\documentstyle[12pt]{amsart}
\newtheorem{Theorem}{Theorem}[section] 
\newtheorem{Lemma}[Theorem]{Lemma}
\newtheorem{Corollary}[Theorem]{Corollary}

\newtheorem{Example}[Theorem]{Example}

\def\char{\operatorname{char}} 
\def\ini{\operatorname{in}} 
\def\gin{\operatorname{gin}}

\def\sk{\smallskip}
\def\q{{\frak q}} 
\def\B{{\cal B}}
\def\a{{\alpha}}

\begin{document}

\title{Borel-fixed ideals and reduction number}

\author{L\^e Tu\^an Hoa and Ng\^o Vi\^et Trung}
\address{Institute of Mathematics,  Box 631, B\`o H\^o, 10000 Hanoi, Vietnam}
\email{lthoa@@thevinh.ncst.ac.vn}
\email{nvtrung@@thevinh.ncst.ac.vn}
\thanks{The authors are partially supported by the National Basic Research Program of Vietnam} 

\maketitle

\section*{Introduction} \sk

Let $A$ be a standard graded algebra over an infinite field $k$. An ideal $\q = (z_1,\ldots,z_s)$, where $z_1,\ldots,z_s$ are linear forms of $A$, is called an {\it $s$-reduction} of $A$ if $\q_t = A_t$ for $t$ large enough. The {\it reduction number} of $A$ with respect to $\q$, written as $r_\q(A)$, is the minimum number $r$ such that $\q_{r+1} = A_{r+1}$. The {\it $s$-reduction number} of $A$ is defined as 
$$ r_s(A) := \min\{r_\q(A)|\ \text{$\q = (z_1,\ldots,z_s)$ is a reduction of
$A$}\}.$$\par
Let $d = \dim A$. It is well-known that a reduction $\q$ of $A$  is  minimal with respect to inclusion if and only if $\q$ can be generated by $d$ elements.  In this case,  $k[z_1,\ldots,z_d] \hookrightarrow A$ is a Noether normalization of $A$ and the reduction number $r_\q(A)$  is the maximum degree of the generators of $A$ as a graded $k[z_1,\ldots,z_d]$-module [V1].  For short, we set $r(A) = r_d(A)$.
The reduction number $r(A)$ can be used as a measure for the complexity of $A$. 
For instance, we can relate $r(A)$ to other important invariants of $A$ such that the degree, the arithmetic degree and the Castelnuovo-Mumford regularity (see [T1], [V1], [V2]). 

Let $I$ be an arbitrary homogeneous ideal in a polynomial ring $R = k[x_1,\ldots,x_n]$. It is shown recently in [C] and [T3] (see also [BH]) that $r(R/I) \le r(R/\ini(I))$, where $\ini(I)$ denotes the initial ideal of $I$ with respect to a given term order.
In particular, we have $r(R/I) = r(R/\gin(I))$, where $\gin(I)$ denotes the generic initial ideal of $I$ with respect to the reverse lexicographic term order [T2]. Since generic initial ideals are Borel-fixed (see the definition in Section 1), we may restrict the study on the reduction number to that of Borel-fixed ideals. 
If $\char(k) = 0$,  Borel-fixed ideals are characterized by the so-called strong stability which gives information on their monomials [BaS]. Similar characterizations can be established for the positive characteristic cases [P]. But these characterizations are not good enough for certain problems. For instance, Conca [C] has raised the question whether $r(R/I) \le r(R/I^{lex})$, where $I^{lex}$ denotes the unique lex-segment ideal  whose Hilbert function is equal to that of $I$. He solved this question for $\char(k) = 0$ by using the strong stability, but his proof does not work for the positive characteristic cases. \par

The aim of this paper is to study the relationship between the $s$-reduction number and Borel-fixed ideals in all characteristics. By definition, Borel-fixed ideals are closed under certain specializations which is similar to the strong stability. Using this property we show that the reduction numbers of $s$-reductions of the quotient ring of a Borel-fixed ideal are attained by $s$-reductions generated by variables (Theorem \ref{reduction}). This gives a pratical way to compute the $s$-reduction number. We will also estimate the number of monomials which can be specialized to a given monomial in the above sense (Theorem \ref{cardinality}). As a consequence, we obtain a combinatorial version of the well-known Eakin-Sathaye's theorem  which estimates the $s$-reduction number by means of the Hilbert function (Corollary \ref{Borel} and Theorem \ref{EaS}). Furthermore, we show 
that the bound of Eakin-Sathaye's theorem is attained by the $s$-reduction 
number when $I$ is a lex-segment monomial ideal (Theorem \ref{Lex}). These results help 
solve Conca's question for all characteristics in a more general setting, namely, that
$r_s(R/I) \le r_s(R/I^{lex})$. Finally, 
since $r(R/I^{lex})$ is extremal in the class of ideals with a given Hilbert function, we will estimate $r(R/I^{lex})$ in terms of some standard invariants of $I$. We shall see that $r(R/I^{lex})$ is bounded exponentially by $r(R/I)$ (Theorem \ref{Lex-bound2}).\par

Thoughout this paper, if $Q \subset R$ is an ideal which generates a reduction of $R/I$, then we will denote its reduction number by $r_Q(R/I)$.

\section{Borel-fixed ideals}

Let $I$ be a monomial ideal of the polynomial ring $R = k[x_1,\ldots,x_n]$.
Let $\cal B$ denote the Borel subgroup of GL$(n,k)$ which consists of the upper triangular invertible matrices. Then $I$ is called a {\it Borel-fixed} ideal if for all $g \in \cal B$, $g(I) = I$. We say that a monomial $x^B$ is a {\it Borel specialization} of a monomial $x^A$ if $x^B$ can be obtained from $x^A$ by replacing every variable $x_i$ of $x^A$ by a variable $x_{j_i}$ with $j_i \le i$. The name comes from the simple fact that any Borel-fixed monomial ideal is closed under Borel specialization.

\begin{Lemma} \label{specialization}
Let $I$ be a Borel-fixed monomial ideal. If $I$ contains $x^A$ then $I$ contains any Borel specialization of $x^A$. 
\end{Lemma}

\begin{pf}
Let $x^B$ be a monomial obtained from $x^A$ by replacing each variable $x_i$ by a variable $x_{j_i}$ with $j_i \le i$, $i = 1,\ldots,n$. Let $g$ be the element of the Borel group $\B$ defined by the linear transformation 
$$ 
g(x_i) = \left\{\begin{array}{lll} x_i & \text{if} & j_i = i,\\  x_i+x_{j_i}  & \text{if} & j_i \neq i. \end{array}\right.
$$
Then $x^B$ is a monomial of $g(x^A)$. Since $g(I) = I$, this implies $x^B \in I$.
\end{pf}

Let $d = \dim R/I$. If $I$ is a Borel-fixed ideal, every associated prime ideals of $I$ has the form $(x_1,\ldots,x_i)$ for $i \ge n- d$ (see e.g. [Ei, Corollary 15.25]). From this it follows that $s$ variables of $R$ generate an $s$-reduction of $R/I$ if and only if they are of the form $x_{i_1},\ldots,x_{i_{s-d}},
x_{n-d+1},\ldots,x_n$ with $1 \le i_1 < \ldots < i_{s-d} \le n- d$. It is clear that $r_{(x_{i_1},\ldots,x_{i_{s-d}},
x_{n-d+1},\ldots,x_n)}(R/I)$ is the least integer $r$ such that all monomials of degree $r+1$ in the remained variables are contained in $I$. 
The following result shows that the computation of the reduction numbers of all $s$-reductions of $R/I$ can be reduced to the above class of $s$-reductions.

\begin{Theorem} \label{reduction}
Let $I$ be a Borel-fixed ideal and $s \ge d = \dim R/I$. Then  \par
{\rm (i)} For every $s$-reduction $\q$ of $R/I$, there exist variables $x_{i_1},\ldots,x_{i_{s-d}}$ with $1 \le i_1 < \ldots < i_{s-d} \le n- d$ such that 
$$r_\q(R/I) = r_{(x_{i_1},\ldots,x_{i_{s-d}},
x_{n-d+1},\ldots,x_n)}(R/I).$$\par
{\rm (ii)} $r_s(R/I) = r_{(x_{n-s+1},\ldots,x_n)}(R/I).$
\end{Theorem}

\begin{pf} Let $y_1,\ldots,y_s$ be linear forms of $R$ which generates $\q$ in $R/I$. Without restriction we may assume that
$$y_i = a_{i1}x_1 + a_{i2}x_2 + \cdots + a_{it_i}x_{t_i}\ (i = 1,\ldots,s)$$
with $a_{it_i} \neq 0$ for different indices $t_1,\ldots,t_s$. Let $g$ be the element of the Borel group $\B$ defined by the linear transformation 
$$ 
g(x_j) = \left\{\begin{array}{lll} x_j & \text{if} & j \not\in \{t_1,\ldots,t_s\},\\
y_i & \text{if} & j = t_i,\ 1 \leq i \leq s.\end{array}\right.
$$
Then $g((x_{t_1},\ldots,x_{t_s})) = g((y_1,\ldots,y_s))$. Since $g(I) = I$, this implies that $x_{t_1},\ldots,x_{t_s}$ generate an $s$-reduction of $R/I$ with
$$r_\q(R/I) = r_{(x_{t_1},\ldots,x_{t_s})}(R/I).$$
As observed before,  $x_{t_1},\ldots,x_{t_s}$ must be of the form $x_{i_1},\ldots,x_{i_{s-d}},
x_{n-d+1},\ldots,x_n$ with $1 \le i_1 < \ldots < i_{s-d} \le n- d$. This proves (i). \par

To prove (ii) choose $\q$ such that $r_s(R/I) = r_\q(R/I)$. By (i) there exist variables $x_{t_1},\ldots,x_{t_s}$ such that  
$r_\q(R/I) = r_{(x_{t_1},\ldots,x_{t_s})}(R/I).$
Note that $r_{(x_{t_1},\ldots,x_{t_s})}(R/I)$ is the least integer $r$ such that 
all monomials of degree $r+1$ in the remaining variables are contained in $I$ and that all monomials of degree $r+1$ in 
$x_1,\ldots,x_{n-s}$ are their Borel specializations. By Lemma \ref{specialization}, the latter monomials are contained in $I$, too. This implies
$$r_{(x_{t_1},\ldots,x_{t_s})}(R/I) \ge r_{(x_{n-s+1},\ldots,x_n)}(R/I) \ge r_s(R/I).$$
So we conclude that $r_s(R/I) = r_{(x_{n-s+1},\ldots,x_n)}(R/I)$. \end{pf}

The case $s = d$ of Theorem \ref{reduction} was already proved by Bresinsky and Hoa [BH, Theorem 11]. They showed that all minimal reductions of $R/I$ have the same reduction number. But their arguments can not be extended to the general case. By Theorem \ref{reduction} (i), there are at most ${n-d \choose s-d}$ different reduction numbers for the $s$-reductions. This number ${n-d \choose s-d}$ can be attained if char$(k) > 0$. This displays a different behaviour than in the case $s = d$.

\begin{Example} \label{reducnumber}
{\rm Assume that $\char(k) = p$. Let $d \leq s < n$ and $1 < a_1<\cdots < a_{n-d}$ be integers. Then 
$$I = (x_1^{p^{a_1}},...,x_{n-s}^{p^{a_{n-d}}}) \subseteq R = k[x_1,...,x_n]$$
is a Borel-fixed ideal. For the $s$-reduction $Q =(x_{i_1},...,x_{i_{s-d}}, x_{n-d+1},...,x_n)$ of $R/I$ with $1 \le i_1 < \ldots < i_{s-d} \le n- d$ we have
$$r_Q(R/I) = p^{a_{j_1}} + \cdots + p^{a_{j_{n-s}}} - n+ s,$$
where $\{ j_1,...,j_{n-s} \} = \{ 1,...,n-d\} \setminus \{i_1,...,i_{s-d}\}$.
Hence the $s$-reductions of $R/I$ have exactly ${n-d \choose s-d}$ different reduction numbers. Moreover, we have
$$r_s(R/I) = p^{a_1} + \cdots + p^{a_{n-s}} - n+ s.$$ }
\end{Example}

If $\char(k) = 0$, Borel-fixed ideals are characterized by a closed property stronger than that of Borel specialization. Recall that a monomial ideal $I$ is called {\it strongly stable} if whenever $x^A \in I$ and $x^A$ is divided by $x_i$, then $x^Ax_j/x_i \in I$ for all $j \le i$.  Any strongly stable monomial ideal is Borel-fixed. The converse holds if $\char(k) = 0$  [BaS, Proposition 2.7]. In this case we can easily compute the reduction number of $R/I$ by the following result. \sk

\begin{Corollary} \label{stability}
Let $I$ be a strongly stable monomial ideal.
For any $s \ge \dim R/I$ we have
$$r_s(R/I) = \min\{t|\  x_{n-s}^{t+1} \in I\}.$$
\end{Corollary}

\begin{pf}
By Theorem \ref{reduction} (ii) we have to prove that
$$r_{(x_{n-s+1},\ldots,x_n)}(R/I) = \min\{t|\  x_{n-s}^{t+1} \in I\}.$$
Hence, it is sufficient to show that if $x_{n-s}^{t+1} \in I$ then all monomials of degree $t+1$ in $x_1,\ldots,x_{n-s}$ are contained in $I$. But this follows from the strong stability of $I$. 
\end{pf}

Example \ref{reducnumber} shows that Lemma \ref{stability} does not hold if $I$ is not strongly stable.\sk

If char$(k)= 0$,  the number of possible reduction numbers for the $s$-reductions of $R/I$ is  much smaller than in the case char$(k) > 0$. In fact, for any $s$-reduction
$Q =(x_{i_1},...,x_{i_{s-d}}, x_{n-d+1},...,x_n)$ with $1 \le i_1 < \cdots < i_{s-d} \le n- d$, we can show similarly as above that
$$r_Q(R/I) = \min\{t|\  x_{j_{n-s}}^{t+1} \in I\},$$
where $j_{n-s}$ is the largest index outside the set $\{i_1,...,i_{s-d},n-d+1,...,n\}$. Since there at most $s-d+1$ such indices, Theorem \ref{reduction} (i) 
shows that there are at most $s-d+1$ different reduction numbers for the $s$-reductions.

\begin{Example} 
{\rm Let $I$ be the ideal generated by all monomials bigger or equal a monomial in the list $x_1^{a_1},...,x_{n-d}^{a_{n-d}}$ with respect to the graded lexicographic order, where $1 < a_1<\cdots < a_{n-d}$. It is easy to see that this ideal is strongly stable and the $s$-reductions of $R/I$ have exactly $s-d+1$ different reduction numbers.}
\end{Example}

The set of all monomials which can be Borel-specialized to $x^A$ will be denoted by $P(x^A)$. If we can estimate the cardinality $|P(x^A)|$ of $P(x^A)$, we can decide when $x^A \in I$, depending on the behavior of the Hilbert function of $I$.

\begin{Lemma} \label{contain}
Let $I$ be a Borel-fixed ideal. Assume that $\dim_k(R/I)_t < |P(x^A)|$ for  $t = \deg x^A$. Then $x^A \in I$.
\end{Lemma}

\begin{pf} 
If $x^A \not\in I$, then $P(x^A) \cap I = \emptyset$ by Lemma \ref{specialization}. Since $P(x^A)$ consists of monomials of degree $t$, this implies $\dim_k(R/I)_t \ge |P(x^A)|$, a contradiction.
\end{pf}

\begin{Theorem}   \label{cardinality}
Suppose $x^A = x_{i_1}^{\a_{i_1}}\cdots x_{i_s}^{\a_{i_s}}$ with $\a_{i_1},\ldots,\a_{i_s} > 0$, $1 \le i_1 < \ldots < i_s \le n$. Put $i_{s+1} = n+1$. Then
$$ |P(x^A)|  \ge \sum_{t=1}^s {\a_{i_1} + \cdots + \a_{i_t}+ i_{t+1}-i_t -1 \choose i_{t+1}-i_t-1}-s+1.$$
\end{Theorem}

\begin{pf} 
The cases $n = 0$ and $\deg x^A = 0$ are trivial because $x^A = 1$. Assume that $n \ge 1$ and $\deg x^A > 0$.\par

If $i_s = n$, we let $x^B =  x_{i_1}^{\a_{i_1}}\cdots x_{i_{s-1}}^{\a_{i_{s-1}}}$ and consider $x^B$ as a monomial in the polynomial ring $S = k[x_1,\ldots,x_{n-1}]$.
Any monomial of $P(x^A)$ is the product of a monomial of $P(x^B) \cap S$ with $x_n^{\a_n}$. The converse also holds. Hence $|P(x^A)| = |P(x^B) \cap S|$. Using induction on $n$ we may assume that
$$|P(x^B) \cap S| \ge \sum_{t=1}^{s-1} {\a_{i_1} + \cdots + \a_{i_t}+i_{t+1}-i_t-1 \choose i_{t+1}-i_t-1}-(s-1)+1.$$
Since $i_{s+1} = n+1 = i_s +1$, we have
$${\a_{i_1} + \cdots + \a_{i_s}+i_{s+1}-i_s -1\choose i_{s+1}-i_s-1} = 1.$$
So we get
$$ |P(x^A)|  = |P(x^B)\cap S| \ge \sum_{t=1}^s {\a_{i_1} + \cdots + \a_{i_t}+i_{t+1}-i_t-1 \choose i_{t+1}-i_t-1}-s+1.$$\par

If $i_s < n$, we divide $P(A)$ into two disjunct parts $P_1$ and $P_2$. The first part $P_1$ consists of monomials divided by $x_{i_1}$,  and the second part $P_2$ consists of monomials not divided by $x_{i_1}$. Set $x^C = x_{i_1}^{\a_1-1}x_{i_2}^{\a_2}\cdots x_{i_s}^{\a_s}$. Every monomial of $P_1$ is the product of $x_{i_1}$ with a monomial of $P(x^C)$. The converse also holds.  Hence $|P_1| = |P(x^C)|$. Using induction on $\deg(x^A)$ we may assume that
\begin{align*}
|P(x^C)| \ge & \sum_{t=1}^s {\a_{i_1} + \cdots + \a_{i_t}+i_{t+1}-i_t-2 \choose i_{t+1}-i_t-1} - s+1\\
\ge &{\a_{i_1} + \cdots + \a_{i_s}+i_{s+1}-i_s -2 \choose i_{s+1}-{i_s}-1} 
\end{align*}
Note that the sum should starts from $t = 2$ to $s$ if $a_{i_1} = 1$.  In this case, the above formula holds because ${\a_{i_1}-2 \choose \a_{i_1}+i_2-i_1-1} = {i_2-i_1-1 \choose 0} = 1$. To estimate $|P_2|$ let $x^D = x_{i_1+1}^{\a_{i_1}}\cdots x_{i_s+1}^{\a_{i_s}}$. It is obvious that every monomial of $P(x^D)$ does not contain $x_{i_1}$ and can be Borel-specialized to $x^A$.
Therefore,  $P(x^D)$ is contained in $P_2$. Using induction on $i_s$ we may assume that
\begin{align*}
|P(x^D)|  \ge & \sum_{t=1}^{s-1} {\a_{i_1} + \cdots + \a_{i_t}+i_{t+1}-i_t-1 \choose i_{t+1}-i_t-1}\\
&  + {\a_{i_1} + \cdots + \a_{i_s}+i_{s+1}-i_s-2\choose
i_{s+1}-i_s-2}- s+1.
\end{align*}
Summing up we obtain
\begin{align*}
|P|  = & |P_1|+|P_2|  \ge  |P(x^C)| + |P(x^D)|\\
\ge &  {\a_{i_1} + \cdots + \a_{i_s}+i_{s+1}-i_s -2 \choose i_{s+1}-i_s-1} +  \sum_{t=1}^{s-1} {\a_{i_1} + \cdots + \a_{i_t}+i_{t+1}-i_t -1 \choose i_{t+1}-i_t-1}\\ & + {\a_{i_1} + \cdots + \a_{i_s}+i_{s+1}-i_s -2 \choose i_{s+1}-i_s-2}-s+1\\
= & \sum_{t=1}^s{\a_{i_1} + \cdots + \a_{i_t}+i_{t+1}-i_t-1 \choose i_{t+1}-i_t-1}-s+1.
\end{align*}
\end{pf}

The bound of Theorem \ref{cardinality} is far from being the best possible as one can realize from the proof. However, it is sharp in many cases.

\begin{Example}
{\rm If $R = k[x_1,x_2,x_3]$ we have
$P(x_1x_3)  = \{x_1x_3,x_2x_3\}.$ Hence}
$$|P(x_1x_3)| =  2 = {3-1+1-1 \choose 1} + {4-3+1-1 \choose 1} -2+1.$$
\end{Example}

An interesting application of Theorem \ref{cardinality} is the following bound for the reduction number.

\begin{Corollary} \label{Borel}
Let $I$ be a Borel-fixed monomial ideal. Assume that
$$\dim_k(R/I)_t < {s+t \choose t}$$
for some integers $s, t \ge 1$. Then $x_{n-s+1},\ldots,x_n$ generates a reduction of $R/I$ with
$$r_{(x_{n-s+1},\ldots,x_n)}(R/I) \le t-1.$$
\end{Corollary}

\begin{pf}
We have to show that the ideal $(I,x_{n-s+1},\ldots,x_n)$ contains every monomial $x^A$ of degree $t$ in $x_1,\ldots,x_{n-s}$. If we write $x^A = x_{i_1}^{\a_{i_1}}\cdots x_{i_s}^{\a_{i_s}}$  with $1 \le i_1 < \ldots < i_s \le n-s$ and $\a_{i_1}+ \cdots + \a_{i_s} = t$, then Theorem \ref{cardinality} gives
$$|P(x^A)|  \ge {n-i_s+t \choose t} \ge {s+t \choose t} > |P(x^A)|.$$
By Lemma \ref{contain}, this implies $x^A \in I$. 
\end{pf}

\section{Eakin-Sathaye's theorem}

Let $R = k[x_1,\ldots,x_n]$ be a polynomial ring over an {\it infinite} field $k$ of arbitrary characteristic. In this section we will deal with the reduction number of $R/I$ for an arbitrary homogeneous ideal $I$. Let us first recall the following theorem of Eakin and Sathaye.

\begin{Theorem} \label{EaS} {\rm [EaS, Theorem 1]}
Let $I$ be an arbitrary homogeneous ideal in $R$. Assume that
$$\dim_k(R/I)_t < {s+t \choose t}$$
for some integers $s, t \ge 1$. Choose $s$ generic linear forms $y_1,\ldots,y_s$, that is in a non-empty open subset of the parameter space of $s$ linear forms of $R$.  Then  $y_1,\ldots,y_s$ generates a reduction of $R/I$ with
$$r_{(y_1,\ldots,y_s)}(R/I) \le t-1.$$
\end{Theorem}

Eakin-Sathaye's theorem provides an efficient way to estimate the reduction number [V2].  We shall see that Corollary \ref{Borel} (though formulated for Borel-fixed ideals and a fixed reduction) is equivalent to Eakin-Sathaye's theorem.   For that we need the following observations.\sk

First, the reduction number of a reduction generated by generic elements is the smallest one among reductions generated by the same number of generators.

\begin{Lemma} \label{generic}
For every integer $s \ge \dim R/I$ choose  
$s$ generic linear forms $y_1,\ldots,y_s$ in $R$. Then
$y_1,\ldots,y_s$ generates a reduction of $R/I$ with
$$r_{(y_1,\ldots,y_s)}(R/I) = r_s(R/I).$$
\end{Lemma}

\begin{pf}
The statement was already proved for the case $s = \dim R$ in [T2, Lemma 4.2].
The proof for arbitrary $s \ge \dim R$ is similar, hence we omit it.
\end{pf}

Secondly, the smallest reduction number does not change when passing to any generic initial ideal.

\begin{Theorem}  \label{initial} 
Let $\gin(I)$ denote the generic initial ideal of $I$ with respect to  the reverse lexicographic term order.  For every integer $s \ge \dim R/I$ we have 
$$r_s(S/I) = r_s(S/\gin(I)).$$
\end{Theorem}

\begin{pf} 
The statement was already proved for the case $s = \dim R$ in [T2, Theorem 4.3].
The case of arbitrary $s \ge \dim R/I$ can be proved in the same manner (though not trivial). 
\end{pf}

Now we are able to show that Eakin-Sathaye's theorem can be deduced from  Corollary \ref{Borel}. Since the proof relies only on properties  of Gr\"obner basis and Borel-fixed ideals,  it can be viewed as a combinatorial proof. \sk

\noindent{\it Combinatorial proof of Theorem \ref{EaS}.} 
By Lemma \ref{generic}, we have to show that
$r_s(R/I) \le t-1$. Let $\gin(I)$ denote the generic initial ideal of $I$ with respect to the reverse lexicographic term order. From the theory of Gr\"obner bases we know that $\gin(I)$ is a Borel-fixed monomial ideal with 
$\dim_k(R/\gin(I))_t = \dim_k(R/I)_t$ (see e.g. [Ei]).
By Corollary \ref{Borel}, the assumption $\dim_k(R/I)_t < {s+t \choose t}$ implies
$$r_s(R/\gin(I)) \le r_{(x_{n-s+1},\ldots,x_n)}(R/\gin(I)) \le t-1.$$
Now, we only need to apply Theorem \ref{initial} to get back to $r_s(R/I)$.
\qed\sk

On the other hand, Corollary \ref{Borel} can be deduced from Eakin-Sathaye's theorem because according to Theorem \ref{reduction} (ii)  and Lemma \ref{generic} we have
$$r_{(x_{n-s+1},\ldots,x_n)}(R/I) = r_s(R/I) = r_{(y_1,\ldots,y_s)}(R/I)$$
for any Borel-fixed ideal $I$. \sk

We shall see that the bound  of Eakin-Sathaye's theorem is attained exactly by lex-segment ideals.
Recall that a {\it lex-segment} ideal is a monomial ideal $I$ such that if $x^A \in I$ then $x^B \in I$ for any monomial $x^B \ge x^A$ with respect to the lexicographic term order. It is easy to see that lex-segment ideals are strongly stable.

\begin{Theorem}\label{Lex}  
Let $I$ be a lex-segment ideal. Then
$$r_s(R/I) = \min\big\{ t| \   \dim_k(R/I)_t < {s+t \choose t} \big\} - 1.$$
\end{Theorem}

\begin{pf} By Theorem \ref{EaS} and Lemma \ref{generic} we have $r_s(R/I) \le r-1$, where
$$r := \min\big\{ t| \   \dim_k(R/I)_t < {s+t \choose t} \big\}.$$ 
It remains to show that $r_s(R/I) \ge r-1$.
Assume to the contrary that $r_s(R/I) < r-1$. By Theorem \ref{reduction} (ii) we have
$r_{(x_{n-s+1},\ldots,x_n)}(R/I) = r_s(R/I) < r-1.$
Using Lemma \ref{stability} we can deduce that $x_{n-s}^{r-1} \in I$. By the definition 
of a lex-segment ideal, this implies that every monomial of degree $r-1$ which involves 
one of the variables $x_1,\ldots,x_{n-s-1}$ is contained in $I$. Equivalently, the monomials of degree $r-1$ not contained in $I$  involve only the $s+1$ variables $x_{n-s},\ldots,x_n$. Since $x_{n-s}^{r-1} \in I$, this implies
$$\dim_k(R/I)_{r-1}   < {s+r-1 \choose r-1}.$$
This contradicts to the definition of $r$. \end{pf}

Given a  homogeneous ideal $I$ in $R$, we denote by $I^{lex}$ the unique lex-segment ideal whose Hilbert function is equal to that of $I$. 
It is well-known that the Betti numbers of $R/I^{lex}$ are extremal in the class of ideals with a given Hilbert function [Bi], [H], [P]. If  char$(k) = 0$, Conca showed
that the reduction number $r(R/I^{lex})$ is extremal in this sense [C, Proposition 10]. He raised the question whether this result holds for all characteristics.  The following result will settle Conca's question in the affirmative.

\begin{Corollary} \label{Conca}
Let $I$ be an arbitrary homogeneous ideal in $R$ and $s \ge \dim R/I$. Then
$$r_s(R/I) \le r_s(R/I^{lex}).$$
\end{Corollary}

\begin{pf} 
According to Theorem \ref{Lex} we have
$$r_s(R/I^{lex}) = \min\big\{ t| \   \dim_k(R/I)_t < {s+t \choose t} \big\} - 1.$$
By Theorem \ref{EaS}, this implies $r_s(R/I) \le r_s(R/I^{lex})$.
\end{pf}

By Corollary \ref{Conca}, $r(R/I^{lex})$ is extremal in the class of ideals with a given Hilbert function. So it is of interest to estimate $r(R/I^{lex})$ in terms of other invariants of $I$.

\begin{Lemma} \label{Lex-bound1} 
Let $I$ be an arbitrary homogeneous ideal in $R$ and $d = \dim R/I \ge 1$. Let $Q$ be an ideal generated by $d$ linear forms of $R$ which forms a reduction in $R/I$.
Put $e = \ell(R/Q+I)$. Then
$$r(R/I^{lex}) \le d(e-2)+1.$$
\end{Lemma}

\begin{pf} 
By [RVV, Theorem 2.2] we know that
$$\dim_k(R/I)_t \le (e-1){t+d-2 \choose d-1} +  {t+d-1 \choose d-1}.$$
For $t = d(e-2)+2$ we have
$$(e-1){de-d  \choose d-1} +  {de-d+1 \choose d-1} <  {de-d+2 \choose d}.$$
Hence the conclusion follows from Theorem  \ref{Lex}.
\end{pf}

We would like to point out that a bound for $r(R/I)$ in terms of $e$ should be smaller. In fact, we always have 
$$r(R/I) \le r_Q(R/I) \le \ell(R/Q+I)-1 = e-1.$$  
If $R/I$ is a Cohen-Macaulay ring, $e$ is equal to the degree (multiplicity) of $I$. If $R/I$ is not a Cohen-Macaulay ring, we may replace $e$ by the extended  (cohomological) degree of $I$ introduced in [DGV].

\begin{Theorem}\label{Lex-bound2} 
Let $I$ be an arbitrary homogeneous ideal in $R$ and $d = \dim R/I \ge 1$. 
Let $a_1 \ge a_2 \ge \cdots \ge a_s$ be the degrees of the minimal homogeneous generators of $I$. Then \par 
{\rm (i)} $r(R/I^{lex}) \leq d\big[\displaystyle {r(R/I)+n-d \choose n-d} - 2\big]+1$, \par
{\rm (ii)} $r(R/I^{lex}) \leq d(a_1   \cdots a_{n-d} -2)+1$.
\end{Theorem}

\begin{pf} Without loss of generality  we may assume that  $Q = (x_{n-d+1},...,x_n)$  forms  a minimal reduction of $R/I$ with $r_Q(R/I)  = r(R/I)$. 
Since $R_t = (Q+I)_t$ for $t \ge r(R/I)+1$, we have
\begin{align*}
\ell(R/Q+I) & \le \sum_{t=0}^{r(R/I)}\dim_k(R/Q+I)_t\\
&  \le \sum_{t=0}^{r(R/I)} \dim_k(R/Q)_t  = {r(R/I)+n-d \choose n-d}.
\end{align*}
Hence (i) follows from Lemma \ref{Lex-bound1}. To prove (ii) we put
$R' = k[x_1,...,x_{n-d}]$ 
and $I' = (I+Q) \cap R'$. Then $I'$ is generated by forms of degrees 
$a'_1 \leq a_1,\  a'_2 \leq a_2 , ...$ and
$\ell(R/Q+I) = \ell(R'/I')$.
By [Bri] we can choose a regular
sequence $f_1,...,f_{n-d}$ in $I'$ such that $\deg (f_i) = a'_i, \  i =1,...,n-d$.
It is well-known that $\ell(R'/(f_1,\ldots,f_{n-d})) = a_1\cdots a_{n-d}$. Hence 
$$\ell(R/Q+I) \le a'_1\cdots a'_{n-d} \le a_1 \cdots a_{n-d}.$$
Thus, (ii) follows from Lemma \ref{Lex-bound1}. 
\end{pf}

Finally we give some examples which show that the bounds of Theorem \ref{Lex-bound2} are sharp.

\begin{Example}
{\rm   Let $I = (x_1,...,x_{n-d})^2$. It is easy to see that $r(R/I) =1$ and 
$$\dim_k(R/I)_t = {d+t-1 \choose d-1} + (n-d){d+t-2 \choose d-1}$$
for all $t \geq 1$. By Theorem \ref{Lex} we have
\begin{align*}
r(R/I^{lex}) & =\min \{ t;\  {d+t-1 \choose d-1} + (n-d){d+t-2 \choose d-1} <
{d+t \choose d}\}-1\\
& = d(n-d-1)+1.
\end{align*}
This is exactly the bound (i) of Theorem \ref{Lex-bound2}.}
\end{Example}

\begin{Example}
{\rm  Consider the one-dimensional ideal $I = (x_1^a) \subset R = k[x_1,x_2],\ a \ge 1.$ We have $\dim_k(R/I)_t = a$ for all $t \ge a-1$. Hence Theorem \ref{Lex} gives
$$r(R/I^{lex}) = \min\{t|\ a < t+1\}-1 = a-1.$$
This shows that the bound (ii) of Theorem \ref{Lex-bound2} is sharp.}
\end{Example}

\section*{References}

\noindent [BaS] D.~Bayer and M.~Stillman, A criterion for detecting $m$-regularity, Invent. Math. 87 (1987), 1-11.\par

\noindent [Bi] A. Bigatti, Upper bounds for the Betti numbers of a given Hilbert function, Comm. Algebra 21 (7) (1993), 2317-2334.\par

\noindent [BrH] H.~Bresinsky and L.~T.~Hoa, On the reduction number of some graded algebras, Proc. Amer. Math. Soc. 127 (1999), 1257-1263.\par

\noindent [Bri] J. Brian\c{c}on,  Sur le degr\'e des relations entre polyn\^omes, 
C. R. Acad. Sci. Paris S\'er. I Math.  297 (1983), 553-556. \par

\noindent [C] A. Conca, Reduction numbers and initial ideals, Proc. Amer. Math. Soc., to appear.\par

\noindent [DGV] L.R.~Doering, T.~Gunston and W.~Vasconcelos, Cohomological
degrees and Hilbert functions of graded modules, Amer. J. Math. 120 (1998), 493--504. \par

\noindent [EaS] P. Eakin and A. Sathaye, Prestable ideals, J. Algebra 41 (1976), 439-454. \par

\noindent [Ei] D.~Eisenbud, Commutative Algebra with a viewpoint toward Algebraic Geometry, Springer, 1994. \par

\noindent [H] H. Hulett, Maximum Betti numbers of homogeneous ideals with a given Hilbert function, Comm. Algebra 21 (7) (1993), 2335-2350.\par

\noindent [NR] D. G. Northcott and D. Rees, Reductions of ideals
in local rings, Proc. Cambridge Philos. Soc. 50 (1954), 145-158.\par

\noindent [P] K. Pardue, Maximal minimal resolutions, Illinois J. Math. 40 (4) (1996), 564-585. \par

\noindent [RVV] M.~E.~Rossi, G.~Valla, and W.~Vasconcelos, Maximal Hilbert
functions, Results in Math. 39 (2001) 99-114.\par

\noindent [T1] N.~V.~Trung, Reduction exponent and degree bound
for the defining equations of graded rings, Proc. Amer. Math.
Soc. 101 (1987), 229-236. \par

\noindent [T2] N.~V.~Trung, Gr\"obner bases, local cohomology and
reduction number, Proc. Amer. Math. Soc. 129 (1) (2001), 9-18.\par

\noindent [T3] N.~V.~Trung, Constructive characterization of the reduction numbers, Compositio Math., to appear. \par

\noindent [V1] W. Vasconcelos, The reduction number of an algebra,
Compositio Math. 106 (1996), 189-197. \par

\noindent [V2] W. Vasconcelos, Reduction numbers of ideals, 
J. Algebra 216 (1999), 652-664. \par

\end{document}